\documentclass[12pt]{amsart}
\usepackage[margin=1in]{geometry}
    \usepackage{amsmath}
    \usepackage{amsfonts}
    \usepackage{graphicx}
    \usepackage{color}
    \usepackage{amssymb}
    \usepackage{amscd}
    \usepackage{psfrag}
    \usepackage{latexsym}
    \usepackage[mathscr]{eucal}

    \newcommand{\co}[2]{\colon\!#1\!\to\!#2}

    \newcommand{\Hyp}{\ensuremath{\mathbf{H}}}

    \newcommand{\inj}[2]{\ensuremath{{\rm{inj}}(#1,#2)}}

    \theoremstyle{plain}
    \newtheorem*{theorem}{Theorem}
    \newtheorem{proposition}[equation]{Proposition}
    \newtheorem{lemma}[equation]{Lemma}

\begin{document}

\frenchspacing

\title{Bounding surface actions on hyperbolic spaces}

\author{Josh Barnard}

\address{Dept.\ of Mathematics \& Statistics\\
	University of South Alabama\\
	Mobile, AL 36688}
\email{jbarnard@southalabama.edu}

\date{\today}

\subjclass[2000]{20F65}

\begin{abstract}
Given an isometric action of the fundamental group of a closed orientable surface on a $\delta$-hyperbolic space, we find a standard generating set whose translation distances are bounded above in terms of the hyperbolicity constant $\delta$, the genus of the surface, and the injectivity radius of the action, which we assume to be strictly positive.
\end{abstract}

\maketitle

\section{Introduction}

We prove the following (see Section~\ref{sec:prelim} for definitions):
\begin{theorem}\label{oldthm}
Suppose $S$ is a closed orientable surface of genus $g\geq 2$ and that $\pi_1(S)$ acts isometrically on an ultracomplete geodesic space $X$ with $\delta$-thin triangles, where the injectivity radius of the action is bounded below by $\epsilon>0$. There is a constant $D$, depending only on $\delta$, $g$, and $\epsilon$, so that for some standard generating set $\{\alpha_1,\ldots,\alpha_{2g}\}$ for $\pi_1(S)$ and some point $x\in X$, we have for $i=1,\ldots,2g$ that the distance from $x$ to $\alpha_i(x)$ is at most $D$.
\end{theorem}

We originally proved a weaker version of this theorem in the appendix to~\cite{Thesis}. The theorem above is used by Bowditch in \cite[Theorem 3.3]{Bow}, where he shows that there are only finitely many atoroidal surface-by-surface bundles with base and fiber of bounded genus. For such a bundle the fundamental group of the base admits an action on the curve complex of the fiber satisfying the hypotheses of Theorem~\ref{oldthm}.

In case $X=\Hyp^n$, a stronger result can be obtained, where one is allowed to take $x$ to be any point lying on any loxodromic axis in $\Hyp^n$. Bonahon~\cite[Lemma 1.10]{Bon} gives an elegantly simple proof of this fact, combining upper bounds on the area of hyperbolic triangles with lower bounds on the injectivity radius of the action to produce the desired upper bounds. In particular, this argument fails for $\delta$-hyperbolic spaces because of the lack of area bounds for geodesic triangles.

A corollary of this theorem is the fact, originally claimed by Gromov~\cite{Gromov} and proven by Sela~\cite{Sela}, that in any hyperbolic group there are only finitely many conjugacy classes of subgroup isomorphic to a given surface group.

More general versions of Theorem~\ref{oldthm} have recently appeared. In~\cite{Bow2} Bowditch shows how to adapt the proof here to apply also to surfaces with boundary, as well as to nonorientable surfaces. Related results using different methods have also been found by Dahmani and Fujiwara~\cite{df}, where $\pi_1(S)$ is replaced by any one-ended subgroup of the mapping class group of $S$ and $X$ is the curve complex of $S$.

\section{Preliminaries}\label{sec:prelim}

A $\delta$-hyperbolic metric space $(X,d_X)$ is said to be \emph{ultracomplete} if there is a geodesic joining any pair of distinct points in either $X$ or its boundary at infinity. If $X$ is proper, then it is ultracomplete. The converse does not hold, as exhibited by curve complexes, for example (see~\cite[Lemma 5.14]{minsky}).

Suppose $S$ is a closed orientable surface of genus $g\geq 2$, and that $\pi_1(S)$ acts on $X$ by isometries. The \emph{injectivity radius} of the action is defined to be
$$
\inj{X}{\pi_1(S)}=\inf_{1\neq h\in\pi_1(S)}\left(\lim_{n\to\infty}\,\frac{d_X(x,h^nx)}{n}\right).
$$
This quantity does not depend on the choice of $x\in X$.

Let $p\co{\widetilde{S}}{S}$ denote the universal cover of $S$, and for any subset $Z\subseteq S$ set $\widetilde{Z}=p^{-1}(Z)$. A \emph{carrier graph} for the action of $\pi_1(S)$ on $X$ is a pair $(\Gamma,\gamma)$ where:
\begin{enumerate}
\item $\Gamma$ is a graph embedded on $S$, the inclusion of which into $S$ induces a surjection $\pi_1(\Gamma)\to\pi_1(S)$;
\item $\gamma\co{\widetilde{\Gamma}}{X}$ is a $\pi_1(S)$-equivariant map sending edges of $\widetilde{\Gamma}$ to rectifiable paths.
\end{enumerate}

From the path-metric $d_{\mathrm{path}}$ induced by $d_X$ on $\gamma(\widetilde{\Gamma})$ we obtain a pseudo-metric $d_{\gamma}$ on $\Gamma$ by declaring for all $a,b\in\Gamma$ that
$$
d_{\gamma}(a,b)=\inf\{ d_{\mathrm{path}}(A,B)\mid A\in\gamma(\widetilde{a}),\,B\in\gamma(\widetilde{b})\}.
$$
Each (open) edge $e\subset\Gamma$ inherits a length $\ell_{\gamma}(e)$ from this metric in the usual way, which we term the $\gamma$-length of $e$, and we set
$$
\ell_{\gamma}(\Gamma)=\sum_{e\subset\Gamma}\ell_{\gamma}(e).
$$

Let $\mathcal{C}$ denote the collection of all carrier graphs for the action of $\pi_1(S)$ on $X$. We say that a carrier graph $(\Gamma,\gamma)$ is \emph{$1$--minimal} if
$$
\ell_{\gamma}(\Gamma)<\inf_{(\Gamma',\gamma')\in\mathcal{C}}\{\ell_{\gamma'}(\Gamma')\}+1.
$$

\section{Main Proposition}\label{sec:prop}

The heavy lifting behind the proof of Theorem~\ref{oldthm} is the following, to whose proof the rest of this section is devoted.
\begin{proposition}
Given $S$ and $X$ as in Theorem~\ref{oldthm}, suppose $(\Gamma,\gamma)$ is a 1--minimal carrier graph for the action of $\pi_1(S)$ on $X$. Then each edge of $\Gamma$ has $\gamma$-length bounded above in terms of the injectivity radius lower bound $\epsilon$, the hyperbolicity constant $\delta$, and the number of edges in $\Gamma$.
\end{proposition}

To prove this proposition, we suppose first that $(\Gamma,\gamma)$ is 1--minimal and has the additional property that $\gamma$ sends each edge of $\widetilde{\Gamma}$ to a geodesic arc in $X$. If we prove the proposition for such carrier graphs, finding that the upper bound is $L$, say, then the result holds for 1--minimal carrier graphs in general with upper bound $L+1$. We may thus assume henceforth for each edge $e$ of $\widetilde{\Gamma}$ that $\gamma(e)$ is a geodesic arc in $X$.

Fix an edge $e$ of $\widetilde{\Gamma}$ and let $n$ denote the number of edges of $\Gamma$. Because $\Gamma$ is embedded $\pi_1$-surjectively, the components of $S\setminus\Gamma$ are simply-connected. Thus we may find a null-homotopic cycle in $\Gamma$, innermost on $S$ and containing $p(e)$. Because $S$ is a surface, such a cycle traverses a given edge of $\Gamma$ no more than twice, and so has combinatorial length no more than $2n$. Let $\Delta$ denote a lift of this cycle to $\widetilde{\Gamma}$ containing $e$.

The key to this proposition is the following lemma.
\begin{lemma}
If $\ell_{X}(\gamma(e))$ is strictly greater than $2(2n-2)\delta+2$, then for any point $s$ on $e$ sufficiently far from its endpoints, there is a homotopically essential loop $\sigma$ on $S$ with $\sigma\cap\Gamma=p(s)$ and having the property that $\gamma$ may be $\pi_1(S)$-equivariantly extended to $\widetilde{\sigma}$ so that $\ell_{\gamma}(\sigma)\leq 2(2n-2)\delta+2$.
\end{lemma}

\begin{proof}
Denote the endpoints of $e$ as $r$ and $u$. For any two points $a,b\in e$, we let $[a,b]$ denote the subsegment of $e$ bounded by $a$ and $b$. In particular we will allow $[r,u]$ to denote $e$ itself.

Let $s$ be any point on $e$ with the property that $\ell_{X}(\gamma([r,s]))$ and $\ell_{X}(\gamma([s,u]))$ are both strictly greater than $(2n-2)\delta+1$. Because $\gamma(\Delta)$ is a geodesic polygon with no more than $2n$ sides, there is a point $s'$ on some other edge $e'$ of $\Delta$ so that $d_X(\gamma(s),\gamma(s'))\leq (2n-2)\delta$. We show first that $p(e)=p(e')$, so that $s$ and $s'$ project to $S$ so as to lie on the same edge of $\Gamma$.

Suppose for contradiction then that $p(e)\neq p(e')$, and let $\alpha$ be a path in $\widetilde{S}$ joining $s$ to $s'$ and otherwise disjoint from $\widetilde{\Gamma}$. We may extend $\gamma$ to $\alpha$ such that $\ell_{X}(\gamma(\alpha))\leq (2n-2)\delta$.

We claim that one of the components of $p(e\setminus s)$ may be removed from $\Gamma\cup p(\alpha)$ so that the remaining graph is still $\pi_1$-surjectively embedded on $S$. To see this, let $e_1$ and $e_2$ denote the components of $e\setminus s$, and suppose $(\Gamma\setminus p(e_i))\cup p(\alpha)$ does not carry all of $\pi_1(S)$ for some $i$. Then there is some essential curve on $S$ intersecting $\Gamma$ only on $p(e_i)$. Thus $e$ has some translate $f$ in $\Delta$, and $e_i$ may be joined to $f$ by a path in $\widetilde{S}$ otherwise disjoint from $\widetilde{\Gamma}$ and $\alpha$. The assumption that $p(e)\neq p(e')$ implies that $e'\neq f$. Thus it is impossible to join both $e_1$ and $e_2$ to $f$ by such paths, because one of the $e_i$ must be separated in $\Delta$ from $f$ by the endpoints of $\alpha$.

If we set $\Gamma'=(\Gamma\setminus p(e_i))\cup p(\alpha)$, choosing $i$ so that $\Gamma'$ carries all of $\pi_1(S)$, we can then equivariantly extend $\gamma$ to the orbit of $\alpha$, sending $\alpha$ to a geodesic. Letting $\gamma'$ denote the restriction of this extension to $\widetilde{\Gamma'}$, this makes $(\Gamma',\gamma')$ a carrier graph for the action of $\pi_1(S)$ on $X$. Because the portion of $p(e)$ removed in forming $\Gamma'$ is strictly longer than $(2n-2)\delta+1$, we have that
\begin{align*}
\ell_{\gamma'}(\Gamma')& =\ell_{\gamma}(\Gamma)+\ell_{\gamma'}(p(\alpha))-\ell_{\gamma'}(p(e_i))\\
& <\ell_{\gamma}(\Gamma)+(2n-2)\delta-\big((2n-2)\delta+1\big)=\ell_{\gamma}(\Gamma)-1,
\end{align*}
contradicting the 1--minimality of $(\Gamma,\gamma)$.

Thus $p(s)$ and $p(s')$ both lie on $p(e)$ in $\Gamma$, and $p(\alpha)$ joins $p(e)$ to itself in the complement of $\Gamma$.
	\begin{figure}[ht]
	\begin{center}
		\includegraphics[width=4in]{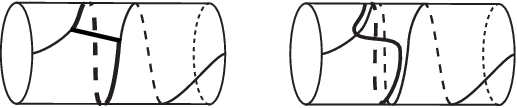}
		\caption{\small{The thin arc is a portion of $\Gamma$, while the bold loops are $p(\alpha)\cup[p(s),p(s')]$ on the left and $\sigma$ on the right.}}
		\label{fig:loop}
	\end{center}
	\end{figure}
Now suppose for contradiction that $\ell_{\gamma}([p(s),p(s')])>(2n-2)\delta+1$. If we let $\Gamma'$ denote the graph obtained from $\Gamma$ by deleting $[p(s),p(s')]$ and inserting $p(\alpha)$, then extending $\gamma$ to the orbit of $\alpha$ by sending $\alpha$ to a geodesic produces a carrier graph $(\Gamma',\gamma')$ with the property that
$$
\ell_{\gamma'}(\Gamma')\leq\ell_{\gamma}(\Gamma)+\ell_{\gamma}(p(\alpha))-\ell_{\gamma}([p(s),p(s')])<\ell_{\gamma}(\Gamma)-1.
$$
This again contradicts the 1--minimality of $(\Gamma,\gamma)$. Thus $\ell_{\gamma}([p(s),p(s')])\leq (2n-2)\delta+1$.

We complete the proof of the lemma by defining $\sigma$ to be the concatenation of $p(\alpha)$ with the segment $[p(s),p(s')]$, pushed off $(p(s),p(s')]$ by an arbitrarily small homotopy (which increases length by less than one). Then $\sigma$ is homotopically nontrivial on $S$ because it does not lift to a loop in $\widetilde{S}$, and $\ell_{\gamma}(\sigma)\leq 2(2n-2)\delta+2$.
\end{proof}

With this lemma in hand, we now address the claim of the proposition, supposing for contradiction that the edge $e$ has the property that
$$
\ell_{\gamma}(p(e))>\big(2(2n-2)\delta+2\big)\big(3\lceil 10\delta/\epsilon\rceil+9/2\big)+9\delta+1.
$$
Then we may choose points $s,t\in e$ ordered as $\{r,s,t,u\}$ on $e$ so that

\vspace*{\jot}
$\ell_{\gamma}(p([r,s]))>(2n-2)\delta+1$,

\vspace*{\jot}
$\ell_{\gamma}(p([s,t]))>\big(2(2n-2)\delta+2\big)\big(3\lceil 10\delta/\epsilon\rceil+7/2\big)+9\delta+1$, and

\vspace*{\jot}
$\ell_{\gamma}(p([t,u]))>(2n-2)\delta+1$.

\vspace*{\jot}\noindent We use the lemma above to find two essential loops $\sigma$ and $\tau$ on $S$ intersecting $\Gamma$ only in $p(s)$ and $p(t)$, respectively, each with $\gamma$-length bounded above by $2(2n-2)\delta+2$. These loops are homotopic on $S$ by construction (they have lifts to $\widetilde{S}$ joining the same two edges of $\Delta$). Each is embedded on $S$, and they are disjoint, because they each have $\gamma$-length bounded above by $2(2n-2)\delta+2$, yet pass through points $p(s)$ and $p(t)$ greater than twice that distance apart.
	\begin{figure}[ht]
	\begin{center}
		\includegraphics[width=4in]{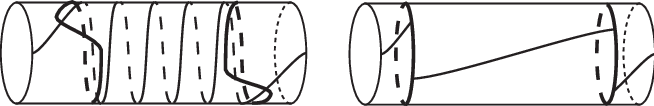}
		\caption{\small{The figure on the left shows a portion of $\Gamma$ on $S$ along with the loops $\sigma$ and $\tau$ in bold. The arc of $\Gamma$ joining them is $c$. The figure on the right shows the corresponding portion of $\Gamma'$, obtained by eliminating $c$, inserting $\sigma$ and $\tau$ (homotoped here for simplicity), and inserting the arc $c'$ cutting across the annulus $A$.}}
		\label{fig:annulus}
	\end{center}
	\end{figure}
It follows that $\sigma$ and $\tau$ cobound an annulus $A$ on $S$. This annulus intersects $\Gamma$ in the arc $c=p([s,t])\subset p(e)$. Let $\hat{\sigma}$ and $\hat{\tau}$ be the components of $\widetilde{\sigma}$ and $\widetilde{\tau}$ containing $s$ and $t$, respectively. We appeal to the following lemma, whose proof we delay, as is it a general result about hyperbolic spaces.
\begin{lemma}
Suppose $h$ is a hyperbolic isometry of the ultracomplete $\delta$-hyperbolic space $X$ with stable translation distance bounded below by $\epsilon>0$, and suppose for some $K\in\mathbf{R}$ that $d_X(x,hx)\leq K$ and $d_X(y,hy)\leq K$ for some points $x,y\in X$. Then there is some $k\in\mathbf{Z}$ so that $d_X(x,h^ky)\leq K\big(3\lceil 10\delta/\epsilon\rceil+7/2\big)+9\delta$.
\end{lemma}
To apply this lemma we take $x=\gamma(s)$ and $y=\gamma(t)$, and we set $K=2(2n-2)\delta+2$, with $h$ denoting the isometry of $X$ represented by $\sigma$ and $\tau$. The lemma now implies that there is some $\hat{t}\in\widetilde{p(t)}\cap\hat{\tau}$ with
$$
d_{X}(\gamma(s),\gamma(\hat{t}))\leq \big(2(2n-2)\delta+2\big)\big(3\lceil 10\delta/\epsilon\rceil+7/2\big)+9\delta.
$$
Therefore there is some arc $c'$ in $A$ joining $p(s)$ and $p(t)=p(\hat{t})$ such that if $\Gamma'=(\Gamma\setminus c)\cup c'$, then $\gamma|_{\Gamma\setminus c}$ extends to a carrier graph $\gamma'\co{\widetilde{\Gamma'}}{X}$ with $\gamma'(\widetilde{c'})$ equal to the $\pi_1(S)$-orbit of a geodesic in $X$ joining $\gamma(s)$ and $\gamma(\hat{t})$.

Observe now that
\begin{align*}
\ell_{\gamma'}(\Gamma')& =\ell_{\gamma}(\Gamma)-\ell_{\gamma}(c)+\ell_{\gamma'}(c')\\
& <\ell_{\gamma}(\Gamma)-\big(2(2n-2)\delta+2\big)\big(3\lceil 10\delta/\epsilon\rceil+7/2\big)-9\delta-1\\
&\qquad\qquad+\big(2(2n-2)\delta+2\big)\big(3\lceil 10\delta/\epsilon\rceil+7/2\big)+9\delta\\
& =\ell_{\gamma}(\Gamma)-1.
\end{align*}
This contradicts the assumed 1--minimality of $(\Gamma,\gamma)$. It follows that
$$
\ell_{\gamma}(e)\leq\big(2(2n-2)\delta+2\big)\big(3\lceil 10\delta/\epsilon\rceil+9/2\big)+9\delta+1,
$$
completing the proof of the Proposition.

\section{Proof of Lemma 2}\label{sec:lem}

We have a hyperbolic isometry $h$ of the ultracomplete $\delta$-hyperbolic space $X$ with stable translation distance bounded below by $\epsilon>0$. We also have points $x$ and $y$ in $X$ such that $d_X(x,hx)\leq K$ and $d_X(y,hy)\leq K$ for some $K\in\mathbf{R}$. We claim there is some $k\in\mathbf{Z}$ so that $d_X(x,h^ky)\leq K\big(3\lceil 10\delta/\epsilon\rceil+7/2\big)+9\delta$.

Because $X$ is ultracomplete, there is a biinfinite geodesic $A_h$ in $X$ that is coarsely preserved by $h$, in the sense that the Hausdorff distance between $A_h$ and $hA_h$ is bounded, so all the translates $h^mA_h$ lie inside the $\delta$-neighborhood of $A_h$. Choose $a\in A_h$ so that $d_X(x,a)=d_X(x,A_h)$. Because the stable translation distance is bounded below by $\epsilon$, if we let $m=\lceil\frac{10\delta}{\epsilon}\rceil+1$, then a geodesic segment $[a,h^ma]$ joining the two indicated points contains a point $b$ with $d_X(b,a)>5\delta$ and $d_X(b,h^ma)>5\delta$.

\bigskip

\noindent{\bf Claim 1}: $d_X(b,A_h)\leq\delta$.

\smallskip

\begin{figure}[ht]
	\begin{center}
		\psfrag{a}{$a$}
		\psfrag{b}{$h^ma$}
		\psfrag{A}{$A_h$}
		\psfrag{B}{$h^mA_h$}
		\psfrag{x}{$x$}
		\psfrag{y}{$h^mx$}
		\psfrag{p}{$b$}
		\includegraphics[width=3.5in]{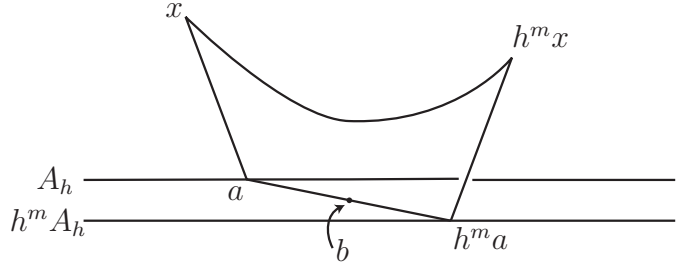}
		\caption{\small{The separating point of Claim 1; the quadrilateral of Claim 2.}}
		\label{fig:claims}
	\end{center}
	\end{figure}

\noindent\emph{Proof of Claim 1:} There are two geodesic triangles sharing the segment $[a,h^ma]$, each of whose other two sides are portions of $A_h$ and $h^mA_h$ sharing an ideal vertex. In each of these triangles there is some point not on $[a,h^ma]$ lying within a distance $\delta$ of $b$. Call these points $b'$ and $b''$, and suppose for contradiction that both lie on $h^mA_h$. Then $b'$ and $b''$ must be separated by $h^ma$ on $h^mA_h$ (see Figure~\ref{fig:claims}). Thus we have
\begin{align*}
2\delta& \geq d_X(b',b)+d_X(b,b'')\geq d_X(b',b'')\\
& =d_X(b',h^ma)+d_X(h^ma,b'')\\
& \geq d_X(h^ma,b)-d_X(b,b')+d_X(h^ma,b)-d_X(b,b'')\geq 8\delta.
\end{align*}
This contradiction establishes Claim 1.

\medskip

\noindent{\bf Claim 2}: $d_X(b,[x,h^mx])\leq 2\delta$.

\smallskip

\noindent\emph{Proof of Claim 2:} The points $a$, $x$, $h^mx$, and $h^ma$ are the vertices of a geodesic quadrilateral in $X$, one of whose sides contains $b$. Thus there is a point $c$ on some other side of this quadrilateral so that $d_X(b,c)\leq 2\delta$. Suppose for contradiction that $c\in[a,x]$. Then
\begin{align*}
d_X(x,c)+3\delta& <d_X(x,c)+d_X(b,a)-d_X(c,b)\\
& \leq d_X(x,c)+d_X(c,a)=d_X(x,a)\\
& =d_X(x,A_h)\\
& \leq d_X(x,c)+d_X(c,b)+d_X(b,A_h)\\
& \leq d_X(x,c)+2\delta+\delta.
\end{align*}
A similar contradiction follows if $c\in[h^ma,h^mx]$. Thus $c\in[x,h^mx]$, establishing Claim 2.

\medskip

Using the point $c\in[x,h^mx]$ from the proof of Claim 2, so that $d_X(b,c)\leq 2\delta$, and choosing $b'$ on $A_h$ with $d_X(b,b')\leq\delta$ as in Claim 1, we have
$$
d_X(x,a)=d_X(x,A_h)\leq d_X(x,c)+d_X(c,b)+d_X(b,b')\leq d_X(x,c)+2\delta+\delta,
$$
from which it follows that
$$
d_X(x,A_h)-3\delta\leq d_X(x,c)\leq d_X(x,h^mx)\leq Km.
$$
In particular, if $d_X(x,hx)\leq K$, then $d_X(x,A_h)\leq Km+3\delta$.

Now for $y$ we have that $d_X(y,hy)\leq K$ implies $d_X(h^ky,h^{k+1}y)\leq K$ for all $k\in\mathbf{Z}$. Applying the argument above, replacing $[x,h^mx]$ with $[h^ky,h^{k+m}y]$, we find that $d_X(h^ky,A_h)\leq Km+3\delta$ for all $k\in\mathbf{Z}$. For each $k$ we let $a_k$ denote a point on $A_h$ nearest to $h^ky$, and note then that
\begin{align*}
d_X(a_k,a_{k+1})& \leq d_X(a_k,h^ky)+d_X(h^ky,h^{k+1}y)+d_X(h^{k+1}y,a_{k+1})\\
& \leq (2m+1)K+6\delta.
\end{align*}
It follows that $\min_{k}d_X(a,a_k)\leq K\big(m+\frac{1}{2}\big)+3\delta$. Thus for the value of $k$ realizing this minimum, we have as claimed that
\begin{align*}
d_X(x,h^ky)& \leq d_X(x,a)+d_X(a,a_k)+d_X(a_k,h^ky)\\
& \leq (Km+3\delta)+\big(K\big(m+\tfrac{1}{2}\big)+3\delta\big)+(Km+3\delta)\\
& \leq K(3m+1/2)+9\delta=K(3\lceil 10\delta/\epsilon\rceil +7/2)+9\delta.
\end{align*}

\section{Proof of Theorem~\ref{oldthm}}\label{sec:thm}

We now show how the main theorem follows from the proposition. Recall that a standard generating set for $\pi_1(S)$ is an ordered set of $2g$ elements $\{\alpha_1,\ldots,\alpha_{2g}\}$ generating $\pi_1(S)$ and satisfying the relation
$$
[\alpha_1,\alpha_2][\alpha_3,\alpha_4]\cdots[\alpha_{2g-1},\alpha_{2g}]=1.
$$
\begin{lemma}
There is a 1--minimal carrier graph $(\Gamma,\gamma)$ for the action of $\pi_1(S)$ on $X$ satisfying the following properties:
\begin{enumerate}
\item $\Gamma$ is a spine on $S$; i.e., the complement $S\setminus\Gamma$ is a single open disk;
\item $\ell_{\gamma}(e)>0$ for all edges $e\in\Gamma$;
\item all vertices of $\Gamma$ are at least trivalent;
\item $\gamma$ sends each edge of $\widetilde{\Gamma}$ to a geodesic arc in $X$.
\end{enumerate}
\end{lemma}

\begin{proof}
Suppose $(\Gamma',\gamma')$ is an arbitrary 1--minimal carrier graph for the action of $\pi_1(S)$ on $X$. We assume that $(\Gamma',\gamma')$ fails each of the conditions above in turn, and show how to produce a carrier graph $(\Gamma,\gamma)$ that satisfies the appropriate condition while remaining 1--minimal.

\smallskip

(i) Because $\pi_1(\Gamma)\to\pi_1(S)$ is surjective, the components of $S\setminus\Gamma$ are topological disks. If this complement has more than one component, then there is at least one edge $e\in\Gamma$ that separates two such components. If $\Gamma'$ is the graph obtained by deleting this edge and $\gamma'$ is the restriction of $\gamma$ to $\Gamma'$, then $(\Gamma',\gamma')$ is a carrier graph for the action of $\pi_1(S)$ on $X$. Moreover $\ell_{\gamma'}(\Gamma')\leq\ell_{\gamma}(\Gamma)$, so that $(\Gamma',\gamma')$ is still 1--minimal.

\smallskip

(ii) Suppose $(\Gamma,\gamma)$ is a $1$--minimal carrier graph with the property that $\ell_{\gamma}(e)=0$ for some edge $e\subset\Gamma$. If $e$ is a loop, then the injectivity radius assumption implies that it must be covered by loops in $\widetilde{\Gamma}$, and hence can be dealt with by the argument of part (i). We therefore assume that $e$ is not a loop, so there is a map $q\co{S}{S}$ which is the composition of the quotient map $S\to S/\{e\}$ and a homeomorphism $S/\{e\}\to S$. Let $\tilde{q}\co{\widetilde{S}}{\widetilde{S}}$ be a map between universal covers that covers $q$. Then there is a carrier graph $(\Gamma',\gamma')$ where $\Gamma'=q(\Gamma)$ and $\gamma'(\tilde{q}(x))=\gamma(x)$ for all $x\in\widetilde{S}$. Moreover it is clear that $\ell_{\gamma'}(\Gamma')\leq\ell_{\gamma}(\Gamma)$, so that $(\Gamma',\gamma')$ is also 1--minimal.

\smallskip

(iii) Any bivalent vertices of $\Gamma$ can be deleted and their adjacent edges merged. If $\Gamma$ has any univalent vertices, the single edges incident to each can be deleted. The resulting carrier graph $(\Gamma',\gamma')$, where $\gamma'$ is the restriction of $\gamma$ to $\widetilde{\Gamma'}$ is clearly 1--minimal if $(\Gamma,\gamma)$ is.

\smallskip

(iv) In this case we take $\Gamma'=\Gamma$, and homotope $\gamma$ so as to map each edge of $\widetilde{\Gamma'}$ to a geodesic arc in $X$ to produce a carrier graph $(\Gamma,\gamma)$ preserving $1$--minimality.
\end{proof}

Suppose $(\Gamma,\gamma)$ is a 1--minimal carrier graph satisfying properties (i)-(iv) above. The proposition implies that each edge $e$ of $\Gamma$ satisfies
$$
\ell_{\gamma}(e)\leq\big(2(2n-2)\delta+2\big)\big(3\lceil 10\delta/\epsilon\rceil+9/2\big)+9\delta+1
$$
where $n$ is the number of edges of $\Gamma$. Using properties (i) and (iii), a simple Euler characteristic argument shows that the number of edges of $\Gamma$ is bounded above by $6g-3$. Thus $\Gamma$ is one of only finitely many graph types, each of which has only finitely many distinct (up to the action of the mapping class group) embeddings into $S$ inducing surjections $\pi_1(\Gamma)\to\pi_1(S)$. In any such graph embedded in $S$ we may find a standard generating set of loops $\{a_1,\ldots,a_{2g}\}$ based at some vertex $v$. By the finiteness properties mentioned above, we may bound the combinatorial length of all such loops by some constant $N_g$ depending only on $g$.

Set $D=N_g\Big(\big(2(12g-8)\delta+2\big)\big(3\lceil 10\delta/\epsilon\rceil+9/2\big)+9\delta+1\Big)$. 
Letting $\alpha_j$ denote the element of $\pi_1(S,v)$ represented by the loop $a_j$, and choosing $x\in\gamma(\tilde{v})$, we have
\begin{align*}
d_X(x,\alpha_j(x))& \leq\ell_{\gamma}(a_j)\\
& \leq N_g\Big(\big(2(12g-8)\delta+2\big)\big(3\lceil 10\delta/\epsilon\rceil+9/2\big)+9\delta+1\big)=D.
\end{align*}

{\bf Acknowledgements}\ \ Thanks are due to Brian Bowditch for many helpful remarks, as well as for the motivation to write this up properly. I am grateful to Daryl Cooper who, among other things, originally suggested the use of carrier graphs, as used in~\cite{Whi}. Many thanks also to the referees for extremely helpful comments.


\begin{thebibliography}{WW}

\bibitem{Thesis}
Barnard, J.: Ends of word-hyperbolic three-manifolds. Ph.D.~thesis, Univ.~of California, Santa Barbara (2004)

\bibitem{Bon}
Bonahon, F.: Bouts des vari\'et\'es hyperboliques de dimension 3. Ann.~Math.~{\bf{124}}, 71--158 (1986)

\bibitem{Bow2}
Bowditch, B.~H.: One-ended subgroups of mapping class groups. Preprint (2007). Available at http://www.warwick.ac.uk/$\sim$masgak/

\bibitem{Bow}
Bowditch, B.~H.: Atoroidal surface bundles over surfaces. GAFA {\bf{19-4}}, 943--988 (2009)

\bibitem{df}
Dahmani, F., Fujiwara, K.: Copies of one-ended groups in mapping class groups. Groups Geom.~Dyn.~{\bf{3}}(3), 359--377 (2009)

\bibitem{Gromov}
Gromov, M.: Hyperbolic groups. In: Essays in Group Theory, {\emph{Math.~Sci.~Res.~Inst.~Publ.}}, vol.~8, pp.~75--263. Springer, New York (1987)

\bibitem{minsky}
Minsky, Y.: The classification of Kleinian surface groups. I. Models and bounds. Ann.~of Math.~(2) {\bf{171}}(1), 1--107 (2010)

\bibitem{Sela}
Sela, Z.: Structure and rigidity in (Gromov) hyperbolic groups and discrete groups in rank 1 Lie groups. II. Geom~Func.~Anal. {\bf{7}}(3), 561--593 (1997)

\bibitem{Whi}
White, M.: Some bounds for closed hyperbolic 3--manifolds. Ph.D.~thesis, Univ.~of California, Santa Barbara (2000)

\end{thebibliography}
\end{document}